\documentclass[final,3p,times]{elsarticle}

\usepackage{amsmath,amsfonts}
\usepackage{graphicx}  
\usepackage{float}  
\usepackage{subfigure}  

\usepackage{amssymb}

\usepackage{soul}
\usepackage{xcolor}
\soulregister{\ref}7
\soulregister{\cite}7

\usepackage{multirow}

\newcommand{\tabincell}[2]{\begin{tabular}{@{}#1@{}}#2\end{tabular}}

\journal{TBD}

\begin{document}

\sethlcolor{yellow}

\begin{sloppypar}

\begin{frontmatter}



\title{Quantifying the Individual Differences of Driver' Risk Perception with Just Four Interpretable Parameters}



\author[label1]{Chen Chen}
\author[label1]{Zhiqian Lan}
\author[label1]{Guojian Zhan}
\author[label1]{Yao Lyu}
\author[label1]{Bingbing Nie}
\author[label1]{Shengbo Eben Li}
\affiliation[label1]{organization={State Key Laboratory of Automotive Safety and Energy, School of Vehicle and Mobility, Tsinghua University},
            city={Beijing},
            country={China}}

\renewcommand{\thefootnote}{}
\footnotetext{Shengbo Eben Li* is the corresponding author. \texttt{Email:lishbo@tsinghua.edu.cn}}
\footnotetext{Codes of this paper can be found here: https://github.com/ChenChenGith/PODAR\_individual\_modeling\_code}

\begin{abstract}
There will be a long time when automated vehicles are mixed with human-driven vehicles. Understanding how drivers assess driving risks and modelling their individual differences are significant for automated vehicles to develop human-like and customized behaviors, so as to gain people's trust and acceptance. However, the reality is that existing driving risk models are developed at a statistical level, and no one scenario-universal driving risk measure can correctly describe risk perception differences among drivers. We proposed a concise yet effective model, called Potential Damage Risk (PODAR) model, which provides a universal and physically meaningful structure for driving risk estimation and is suitable for general non-collision and collision scenes. In this paper, based on an open-accessed dataset collected from an obstacle avoidance experiment, four physical-interpretable parameters in PODAR, including prediction horizon, damage scale, temporal attenuation, and spatial attention, are calibrated and consequently individual risk perception models are established for each driver. The results prove the capacity and potential of PODAR to model individual differences in perceived driving risk, laying the foundation for autonomous driving to develop human-like behaviors.
\end{abstract}




\end{frontmatter}


\section{Introduction}
\label{introduction}
Autonomous driving technology has been developing rapidly in recent years, and some manufacturers are preparing to deploy SAE-L4 on public roads. It is expected to have a mix of automated and regular vehicles in a long future period\cite{ref1}. The interaction performances of automated vehicles, especially on the aspect of driving safety, will be inspected and judged by surrounding human drivers and, consequently, determine their level of trust and cooperation with automated vehicles\cite{ref2}. Research suggests that the more human-like behaviors automated vehicle performs, the more people's trust and acceptance\cite{ref2_1}. Therefore, modelling and knowing how drivers assess driving risks from road objects is essential to support automated vehicles and develop more human-like behaviors, accordingly showing people their higher reliability in terms of safety\cite{ref3}.

Practically, autonomous driving technology is at the vehicle level, and its decision, planning, and control modules collaboratively act as the substitute for a human driver. An proper model that can estimate correct driving risk values is imperative in designing automated vehicle's decision and control functions. However, to the best of our knowledge, a general perceived risk model has not been proposed in previous studies. The parameter most used in current self-driving decision-making and planning algorithms to guarantee safety is a simple measure, i.e., relative distance to surrounding obstacles or traffic participants\cite{ref_dnp_1,ref_dnp_2,ref_dnp_3,ref_dnp_4,ref_dnp_5,ref_dnp_6}. Taking distance as cost or constraint in the optimization process can meet the safety requirements of non-collision. However, the absence of a risk model makes it inadequate to generate behaviors at a proper and acceptable risk level in non-conflict scenes.

More importantly, although some researchers have concluded that automated vehicles were required to be safer than human-driving on average\cite{ref_accep_1,ref_accep_2,ref_accep_3}, one study points out that people would use their own perceived driving ability rather than the safer-than-the-average-driver benchmark\cite{ref_accep_4}, to determine their trust on automated vehicles. That means that the driving risk model should have capacity to distinguish the differences between drivers on risk perception, thus fitting their individual features. This is also of significance to the future customized autonomous driving services.

Despite many efforts in driving risk analysis, driving risk modelling is typically developed at a statistical level. The most population way is to assess risks via analyzing driving behaviors based on one widely proven fact of the significant contribution of driving behaviors to road accidents\cite{ref_behav_1,ref_behav_2,ref_behav_3}. Using behavior data collected from driving simulator\cite{ref_simu_1,ref_simu_2} or naturalistic driving\cite{ref_nat_drv_2}, factors including age and gender\cite{ref6,ref7,ref8}, personality\cite{ref9}, drug-use\cite{ref_drug_1,ref_drug_2,ref_drug_3,ref_drug_4}, and road geometric features\cite{ref12} are examined to obtain their influence on driving risks. The influencing analysis results can support the countermeasures' development for driving safety improvement, but insufficient for directly helping the skill improvement of the individual driver, let alone supporting the decision and planning technology evolution for automated vehicles.

One forward step in grasping risks differences between drivers is the introduction of driving style, which categorizes drivers' habits into different classes as the representation of their driving risk levels\cite{ref_style_1}. Typically, 3-5 classes, for example, aggressive, moderate, and conservative, are defined and correspond to high, normal, and low collision risk levels\cite{ref_style_3,ref_style_4}. Data for driving style classification mainly comes from naturalistic driving collected by onboard diagnostic (OBD), Global Positioning System (GPS), smartphone accelerometer\cite{ref_style_5,ref_style_6}, etc. Due to the lack of labels, unsupervised methods such as clustering, principal component analysis\cite{ref_style_5}, and fuzzy classification\cite{ref_style_2} are usually used based on behavior features like the standard deviation of the speed, longitudinal and lateral acceleration, etc. Another way is to detect specific risky driving events, such as sharp acceleration, braking, and turning, and count their frequencies to determine the driver's driving style class\cite{ref_style_7,ref_style_8,ref_style_9,ref_style_10}. Nevertheless, driving style is still a qualitative concept and always concerns with self-vehicle's motions while neglecting the interactions with other traffic participants, which is far away from the needs of autonomous driving risk assessment.

Of course, there are some surrogate measures that can assess driving risks. Traditional measures like Time-to Collision (TTC)\cite{ref_ttc}, Deceleration Rate to Avoid a Crash (DRAC)\cite{ref_drac_1,ref_drac_2}, Time Exposed Time-to-collision (TET)\cite{ref_tet}, and Post-Encroachment Time (PET)\cite{ref_pet_1,ref_pet_2} perform well in longitudinal and lateral conflict situations like car-following, lane change, but are not suitable in none-conflict scenarios. Recently proposed scene-generalized models like Driver Risk Field (DRF)\cite{ref3,ref_drf}, Safety Field (SF)\cite{ref_sf_1,ref_sf_2} provide ways to model risks in non-conflict situations. However, they consider the motion of only the host vehicle or only the surrounding vehicle in the conflict parties instead of predicting both of them, which will lead to false results in some lateral conflict scenes\cite{ref_podar}. Besides, there need huge efforts on the mass of parameters calibration and high computational load due to field-based representation, which limits their practical application in self-driving technology development.

To correctly and accurately assess driving risk, a universal driving risk model called Potential Damage Risk (PODAR) was proposed and proved to be more feasible and rational under general traffic scenarios than pre-mentioned measures in our previous works\cite{ref_podar}. The PODAR model was constructed according to a physically interpretable process like human cognition, which regards the driving risk as the projection of future potential collision damage. Just a few key parameters are needed to model different characteristics of drivers in risk assessment. By comparing the parameters among the individual driver' PODAR model, the reason why drivers perceive different levels of risk, and, consequently, why drivers perform different driving behaviors when encountering obstacles can be explained.

In this paper, we establish eight drivers' individual driving risk models based on a dataset published by a previous study\cite{ref3}. Four physical-interpretable parameters are calibrated for each driver, including prediction horizon, damage scale, temporal attenuation, and spatial attention. The results show the capacity of the PODAR model in delineating the differences in drivers' risk perception at a physical level. It also blazes a way to design customized driving risk configuration via modifying a few parameters, supporting the design of decision and planning functions for automated vehicles to perform personalized human-like behaviors.

\section{Material and Method}
\label{method}
This section first introduces the the PODAR model. Then, the driving simulation data resource and calibration method are described.

\subsection{PODAR Model}
Risk is defined as the ``likelihood and severity of hazardous events"\cite{ref_risk_def_1,ref_risk_def_2}. To put it more clearly, we suggested an intuitive description that risk is the threat of possible damage, which implies an assumption that damage is the supremum of risk numerically, and the risk is the projection of damage due to motions' uncertainty. Based on this idea, a general-scene available driving risk model, Potential Damage Risk (PODAR), was proposed and verified in our previous work\cite{ref_podar}. 

\begin{figure}[t]
	\centering
	\begin{minipage}[t]{0.48\textwidth}
		\centering
		\includegraphics[width=1\textwidth]{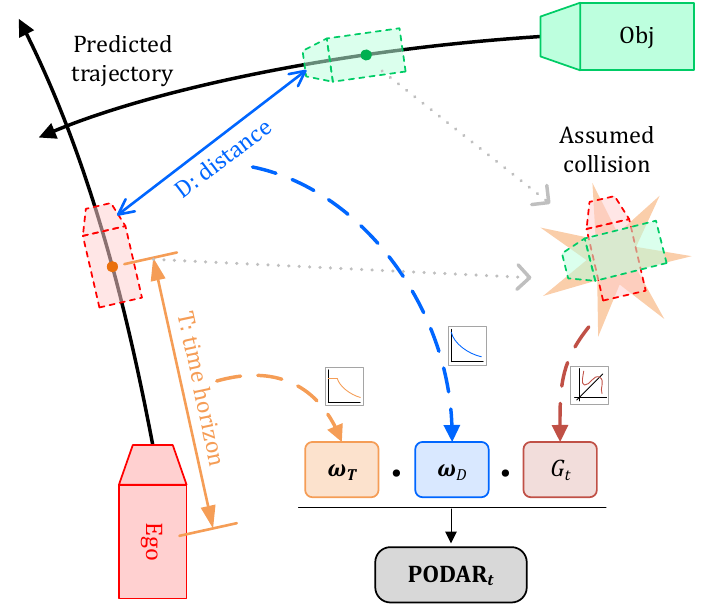}
		\caption{PODAR Model sketch}
		\label{fig_PODAR_Structure}
	\end{minipage}	
	\begin{minipage}[t]{0.48\textwidth}
		\centering
		\includegraphics[width=1\textwidth]{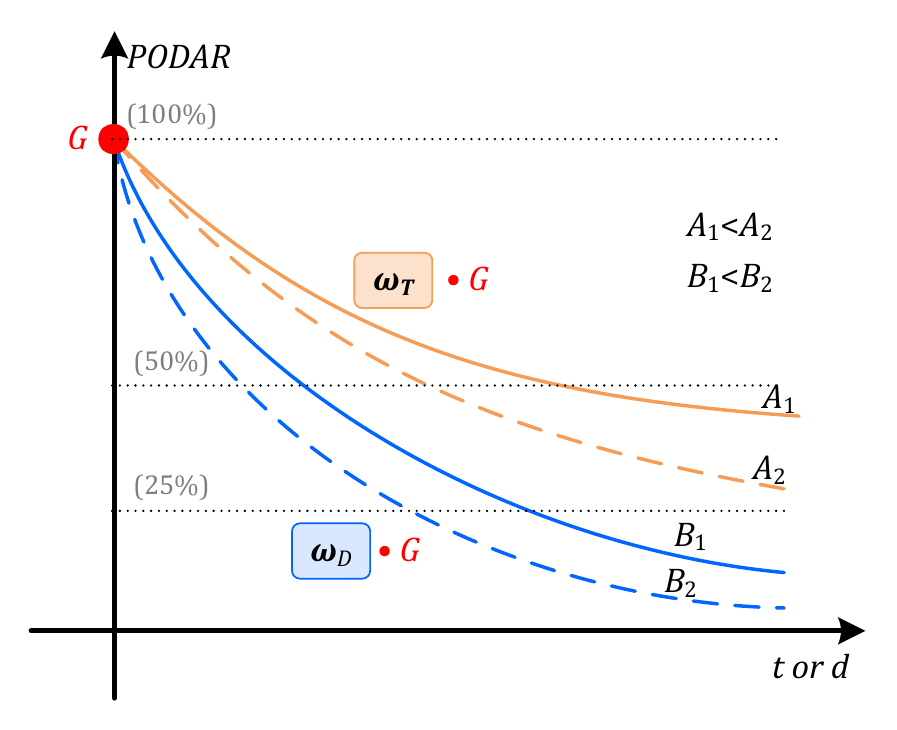}
		\caption{Attenuation curves}
		\label{fig_attenuation_curve}
	\end{minipage}
\end{figure}

The sketch of the PODAR model is shown in Figure \ref{fig_PODAR_Structure}, which implies four key components. 1) First, PODAR predicts the future trajectories of both the host vehicle and surrounding objects. 2) Then, at each predictive timestep, a virtual collision is assumed and thus results in potential damage $G$, serving as the ``severity" in risk definition. 3) Meanwhile, two attenuation functions, $\omega_D$ and $\omega_T$, that respectively describe the damage discount in spatial and temporal dimensionalities are superposed to transfer damage values into perceived risk values, serving as the ``likelihood" in risk definition. 4) Finally, since multiple objects and prediction steps exist, the maximum risk value is determined as the final driving risk to simulate finite human attention. The general structure of the PODAR model is shown as Equation \ref{eq_PODAR_struc}:
\begin{gather}
\label{eq_PODAR_struc}
\text{PODAR}=\max_{t,n}\{\text{PODAR}_{t}^{n}\}\\
\text{PODAR}_{t}^{n}=G_{t}^{n}\cdot\omega_{D}\cdot\omega_{T}\notag\\
t\in\{0,1,...,T\},n\in\{1,2,...,N\}\notag
\end{gather}
where, $\text{PODAR}$ is the finally determined driving risk value at current moment; $t$ is the timestep in prediction horizon $T$; $n$ is the index of surrounding $N$ objects like motor vehicle, pedestrian and bicycle; $\text{PODAR}_{t}^{n}$ is the driving risk from the $n$-th object at timestep $t$; $G_{t}^{n}$ is the potential collision damage from the $n$-th object at timestep $t$; $\omega_{D}$ and $\omega_{T}$ are the spatial and temporal attenuation coefficients related to distance and time, respectively.

Two components in the PODAR model are the potential damage and attenuation functions. Damage $G_{t}^{n}$ comes from an assumed virtual collision with the states of the host vehicle and surrounding object at that predictive timestep, and it can be formulated by Equation \ref{eq_damage}:
\begin{gather}
\label{eq_damage}
G_{t}^{n}=f_g\left(\boldsymbol{v}_{t}^{\text{host}},\boldsymbol{v}_{t}^n,\boldsymbol{p}_{t}^{\text{host}},\boldsymbol{p}_{t}^{n}, M^{\text{host}}, M^n \right) 
\end{gather}
where, $f_g$ is the function for estimating potential collision damage; $\boldsymbol{v}_{t}$ and $\boldsymbol{p}_{t}$ are the velocity vector and the position vector at timestep $t$ of the host vehicle or $n$-th object; $M$ is the virtual mass to emphasize the importance of some vulnerable participants. For practicability, $f_g$ is designed as a kinetic energy-like form, see Equation \ref{eq_damage_exp}, \ref{eq_speed_exp}, and \ref{eq_M_exp}:
\begin{align}
    \label{eq_damage_exp} G_t = & \frac{1}{2}\cdot \left( M^{\text{host}}+M^n \right) \cdot V_t\cdot \left| V_t \right|\\
    \begin{split}
        \label{eq_speed_exp} 
        V_t = & \alpha \cdot \left< \boldsymbol{v}_t^n-\boldsymbol{v}_t^{\text{host}}, \frac {\boldsymbol{p}_t^{\text{host}}-\boldsymbol{p}_t^n}{\left| \boldsymbol{p}_t^{\text{host}}-\boldsymbol{p}_t^n \right|} \right> + \left( 1-\alpha \right) \cdot \left( \left| \boldsymbol{v}_t^{n} \right| + \left| \boldsymbol{v}_t^{\text{host}} \right| \right)
    \end{split} \\
    \label{eq_M_exp} M = & m \cdot s
\end{align}
where, $V_t$ is the integrated representation of host and target object velocities; $\alpha$ is the weight for balancing velocity differential and magnitude; $m$ is the actual mass; and $s$ is the damage sensitivities for different object types. The kinetic energy-like form can give a rational description of all the driving situations and has been verified in our previous study\cite{ref_podar}.

The attenuation functions $\omega_{D}$ and $\omega_{T}$ should be decreasing functions that meet Equation \ref{eq_attenuation}. For practicability, empirical models using negative exponential functions were established as shown in Equation \ref{eq_attenuation_exp1}, \ref{eq_attenuation_exp2}, and Figure \ref{fig_attenuation_curve}.
\begin{align}
	\label{eq_attenuation}
	\omega(\cdot) =&\omega_D (d) \ or \ \omega_T (t) \\
     \text{s.t.} \ \ & \omega(0) =1 \notag\\
    &\omega(x_1)\geq\omega(x_2), \forall{x_1 \leq x_2} \notag\\
	\label{eq_attenuation_exp1} 
	\omega_T =& e^{-A\cdot t}, t\geq0, A\geq0\\
	\label{eq_attenuation_exp2} 
    \omega_D =& e^{-B\cdot d}, d\geq0, B\geq0
\end{align}
where, $d$ is the relative distance between the host vehicle contour and the surrounding object contours. 

$A$ and $B$ are the parameters to be calibrated for each driver. A bigger value will lead to a steeper descent curve, which means that only the near-future and close-distance objects will cause high risks, while far-future collisions and long-distance objects would be neglected. Thus, by analyzing the $A$ and $B$ values difference between drivers, the sensitiveness of drivers to potential damages from the temporal and spatial dimensionalities can be grasped, and the mechanism of the risk perception differences among drivers can be explained.

\subsection{Data Resource}
Data used in this paper comes from the supplementary data of paper \cite{ref3} that proposed the Driver's Risk Field (DRF) model. Here, we present a concise description of  the driving simulation experiment design and data collection.

\begin{figure}[t]
	\centering
	\subfigure[road and obstacle positions design]{
		\includegraphics[width=0.48\linewidth]{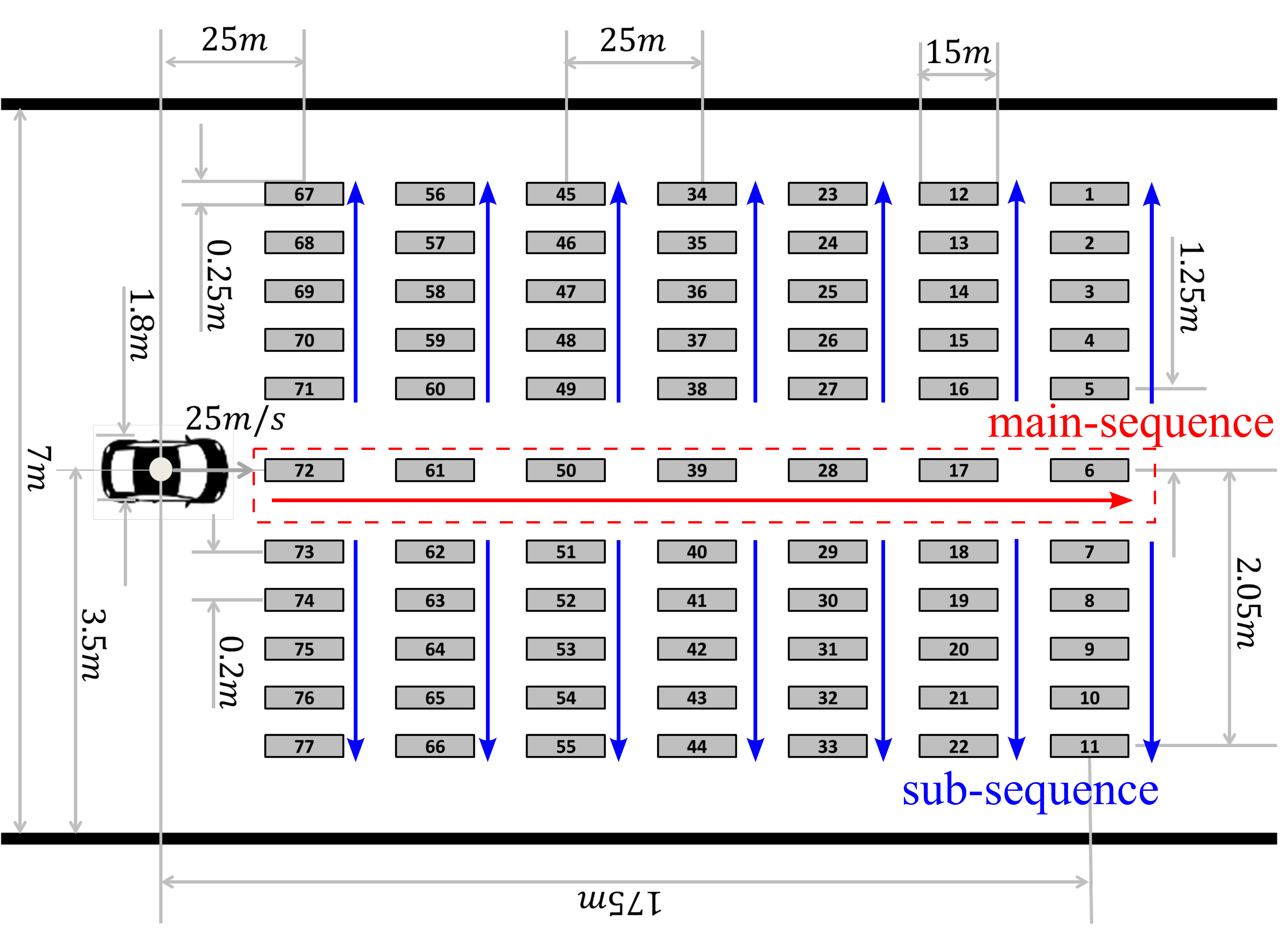}}
	\subfigure[driving task design]{
		\includegraphics[width=0.48\linewidth]{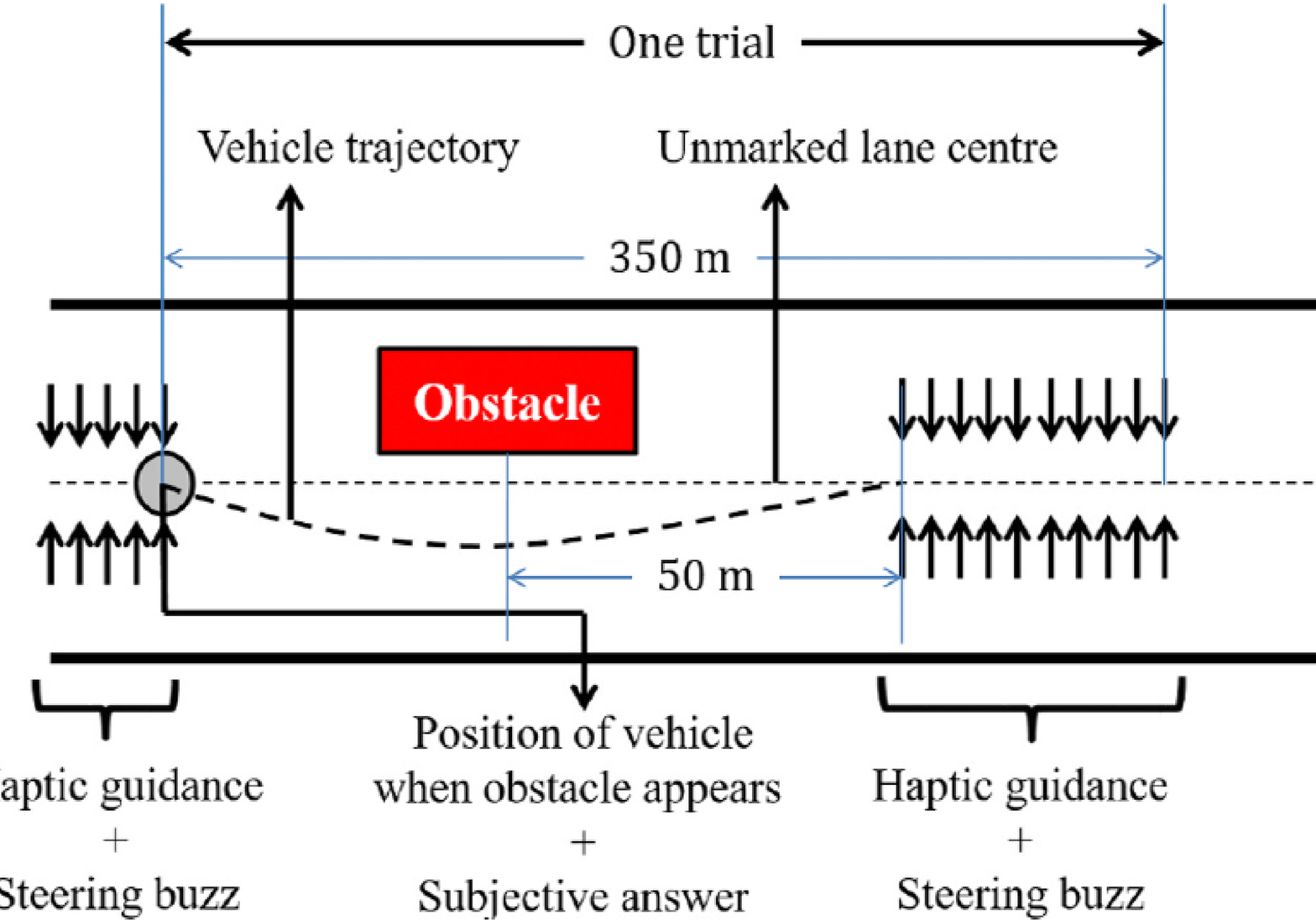}}
	\caption{Driving simulation design(Fig.3 and Fig.4 in literature \cite{ref3})}
	\label{fig_DS_design}
\end{figure}

There were 8 participants (P1-P8) who each conducted 308 (77 obstacles * 4 times) obstacle avoidance trials. The design of the road, obstacle positions and trial section, are shown in Figure \ref{fig_DS_design}. Theoretically, only seven obstacles (ID=O6, O7, O28, O39, O50, O61, and O72, called ``\textbf{main-sequence obstacles}", labeled by red arrow in Figure \ref{fig_DS_design}) located at the center of the road needed to be avoided, and others (called ``\textbf{sub-sequence obstacles}", labeled by blue arrow in Figure \ref{fig_DS_design}) were beyond the width of the host vehicle. The host vehicle drove at a constant speed (25m/s), and guidance torques were exerted on the steering wheel to keep the vehicle staying at the center of the road if it were not in the trail section. When the host vehicle achieved the start of a trial section, one obstacle randomly elected from the 77 obstacles appeared and the guidance torques were deactivated. Drivers were required to manipulate the steering wheel to avoid possible collisions and say aloud a non-negative real number once the obstacle appeared to answer the question that ``How much steering do you think you need, at this moment?". The guidance torques will be re-activated after passing the obstacle for 50m. It should be noted that drivers were instructed to report zero if they thought that no steering wheel operation was needed, and no scale was provided to limit the highest response numbers.

The dataset from the driving simulation experiment could be accessed online. Two signals, the maximum absolute value of the steering angle and oral response numbers, were regarded as the objective and subjective indicators in assessing driver perceived risks. Similar to the treatment in the original literature\cite{ref3}, only the perceived risks at the instant time when obstacles appeared (i.e., the host vehicle achieved the start of the trail section) were focused on. The maximum absolute values of the steering angle were collected during the first 1 second after obstacles appeared. The average of the maximum absolute value of steering angle in four trails, presented as MSA(maximum steering angle), was finally used as the objective signal for each obstacle. The oral response numbers, presented as ORN, had only one record for each obstacle and were not repeated in the following trial. Therefore, a calibration dataset containing 1232 records (77 obstacles * 8 participants * 2 signals) was established and was used to establish the individual driver PODAR models in this paper.

\subsection{Individual PODAR Model Calibration}
\subsubsection*{Data preview and standardization}
\begin{figure}[t]
	\centering
	\subfigure[objective signals responding to 77 obstacles]{
		\includegraphics[width=0.48\linewidth]{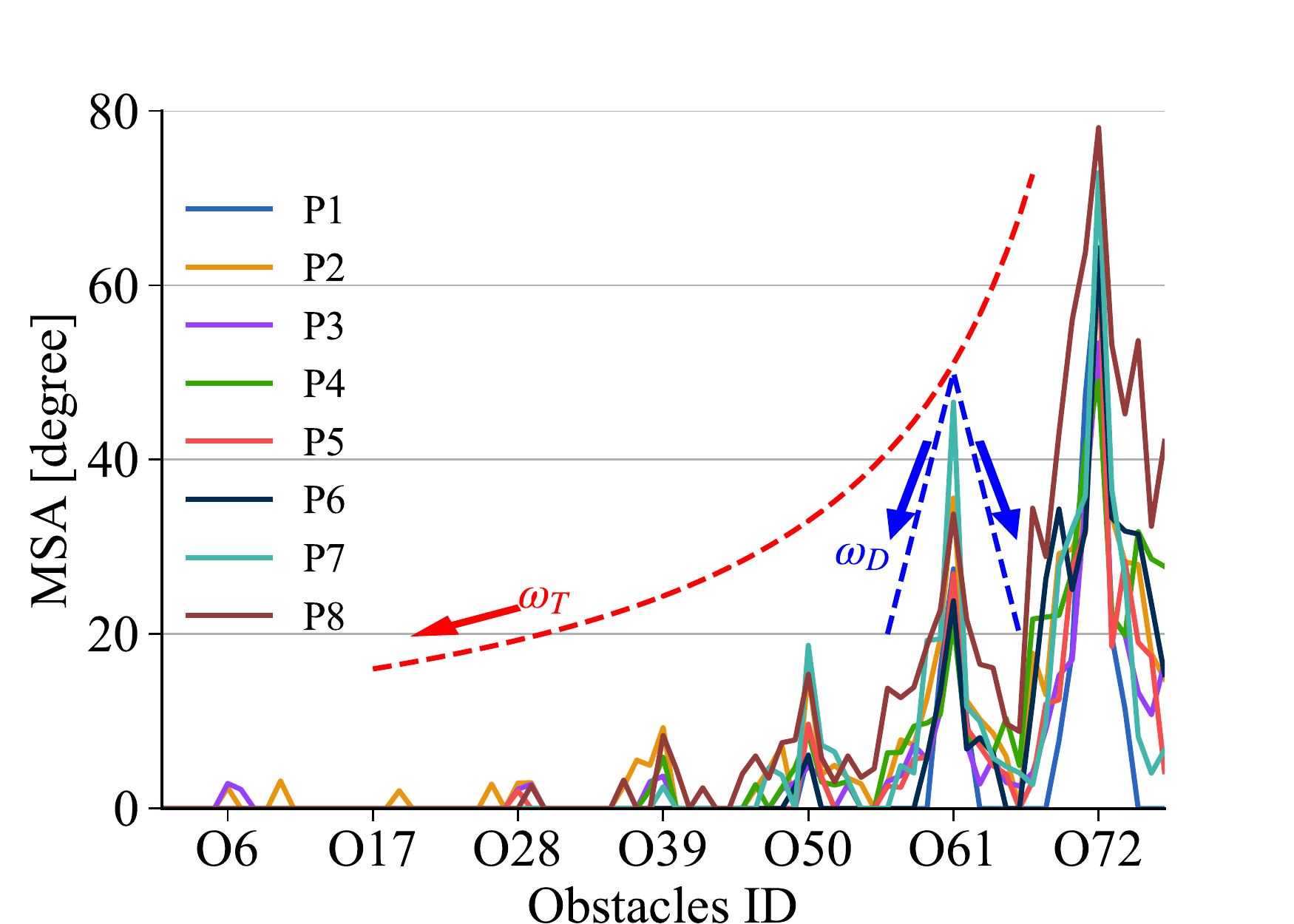}}
	\subfigure[maximum MSAs and ORNs of eight drivers]{
		\includegraphics[width=0.48\linewidth]{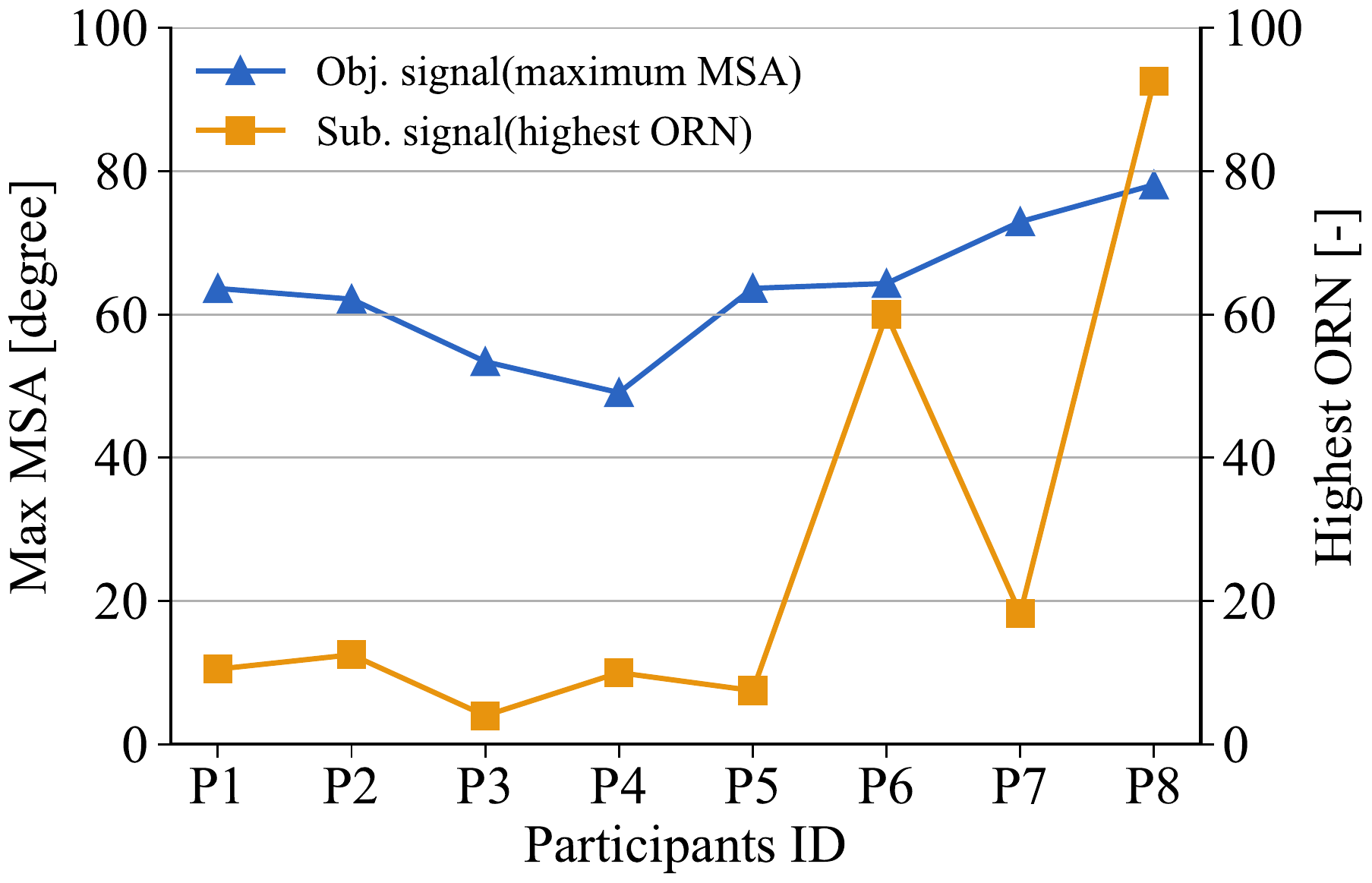}}
	\caption{Calibration dataset overview}
	\label{fig_data_preview}
\end{figure}

Similar to literature \cite{ref3}, MSAs less than 2$^\circ$ were deemed non-conscious steering action and set to 0. Figure \ref{fig_data_preview}a shows the MSAs of 8 drivers to avoid collision with appeared obstacles. Main-sequence obstacles caused peak MSAs, and O72 had the highest value because it appeared at the nearest position, as shown in Figure \ref{fig_data_preview}a. Besides, variations of maximum MSA among drivers were also observed, as the blue line shown in \ref{fig_data_preview}b. MSA is the objective signal, so that it can be used in the comparison between drivers. Thus, standardization was first conducted to normalize the MSA values into the range [0,1], with the maximum MSA among 8 drivers mapping to 1. The subjective signals showed different tendencies because no highest restraint of the oral response numbers was informed to participants, which resulted in large variations, as the yellow line shown in Figure \ref{fig_data_preview}b. Hence, the standardization was conducted on each participant's data instead of the entire dataset.

At the moments when obstacles appeared, the host vehicle was predicted (the first step in the PODAR model) to drive straightly and will have collision ($d_t=0$) with main-sequence obstacles. Therefore, damage attenuations of the main-sequence obstacles will only come from the temporal dimension according to the PODAR structure (Equation \ref{eq_PODAR_struc}). It is easy to imagine that tendencies of MSA or ORN values caused by main-sequence obstacles (red arrow in Figure \ref{fig_data_preview}a) were mainly determined by the parameter $A$ in $\omega_T$, while the tendencies caused by sub-sequence obstacles (blue arrow in Figure \ref{fig_data_preview}a) in the same columns were mainly determined by the parameter $B$ in $\omega_D$. Here, a symmetry assumption was made in the spatial dimension.

\subsubsection*{Parameters needing calibration}
As mentioned above, the key parameters waiting for calibration are the temporal attenuation coefficient $A$, and the spatial attenuation coefficient $B$. 

Besides, there is an implicit parameter, prediction horizon $T$, in the trajectory prediction step of the PODAR model. In our previous works, $T$ was set to 3 seconds and further future motions after that were neglected. It is easy to know that the furthermost prediction horizon is 7 seconds (the farthest obstacles were placed at 175m as shown in Figure \ref{fig_DS_design}a with the host vehicle speed of 25m/s) in the calibration dataset. Whereas, Figure \ref{fig_data_preview}a shows the differences between drivers on the precognition characteristic. We subjectively determined the prediction horizons of each driver for both the objective and subjective signals by examining the minimum number of consecutive MSAs peaks among the main-sequence obstacles. It should be noted that we used a constant speed model to predict the trajectory of the host vehicle because it was expected to drive in a straight line when the obstacle appeared and then drivers controlled the steering according to their perceived risks. 

$M$ of the obstacle for calculating the potential damage $G_{t}^{n}$ in the PODAR model were assigned a fixed value (1.8kg). A kinetic energy form was used to model the $f_g$ as shown in Equation \ref{eq_damage_exp}, \ref{eq_speed_exp}, and \ref{eq_M_exp}, which was the same as that in our previous paper \cite{ref_podar} and the coefficient $\alpha$ was set as 0.7. In fact,the $\alpha$ value would not influence the relative size of estimated damage values of different obstacles by Equation \ref{eq_damage_exp}, because we only concern on the risk at the time when the host vehicle achieves the start of the trial section, where the velocity differential and magnitude at each prediction step in Equation \ref{eq_speed_exp} will be the same. However, this treatment means that the damage estimation is the same for drivers, which is inconsistent with the observation of the maximum MSAs among drivers. Thus, an additional parameter, $k$, was added to scale the damage estimation (Equation \ref{eq_PODAR_struc_adjust}) to show the differences in damage estimation of drivers:
\begin{gather}
\label{eq_PODAR_struc_adjust}
\text{PODAR}_{t}^{n}=(k\cdot G_{t}^{n})\cdot\omega_{D}\cdot\omega_{T}
\end{gather}

The subjective signals can not support the comparison between drivers because that participants were not informed of the scales of ORNs. The parameters $k$ for subjective signal only helped to adjust damage-risk attenuation curves to fit each participant's responses.

To summarize, four parameters needed calibration for each participant: $T$, $k$, $A$, and $B$, corresponding to the prediction horizon, potential damage estimation, temporal and spatial attenuation. The maximum function in Equation \ref{eq_PODAR_struc} for simulating finite human attention was reserved and not considered modifications in the calibration process. Similar to the literature \cite{ref3} $\text{R}^2$ index was calculated to assess how well the PODAR model fits perceived risk, and to compare the performance between the PODAR model and DRF model.

\subsubsection*{Calibration method}
Gradient descent was used to calibrate the parameters for the eight participants one by one. The batch size was 77, i.e., data of all the 77 obstacles was sent into the model for gradient back-propagation in each iteration. The iteration times of the training were set as 100,000. The optimizer was Adam \cite{ref_adam} with an initial learning rate of 0.0001. 

\section{Results}
\label{results}
\subsection{Prediction Horizons}
The prediction horizon ($T$) of each participant are shown in Figure \ref{fig_res_pred_hori}. Almost all drivers had further ahead foresight in subjective than in objective. ORNs were collected once the obstacle appeared, while MSAs were counted according to the manipulations to avoid collisions or reduce perceived risks during 1 second after the obstacle appeared. Differences between the two signals would attribute to the instantaneous and delayed collection methods: drivers would become calm when manipulating the steering while the ORNs were responses to instantaneous stimulation. The shortest and average prediction horizons of objective signals are 3 and 4.2 seconds, respectively, indicating that motions within future 3 seconds are key to the driving risk estimation, and 4 seconds would be a proper horizon to modelling an ``average person" in practice. The subjective prediction horizon is 4 secondes at least and nearly 6 seconds in average.

\subsection{Damage Scale Coefficients}
\begin{figure}[t]
	\centering
	\begin{minipage}[t]{0.48\textwidth}
		\centering
		\includegraphics[width=1\textwidth]{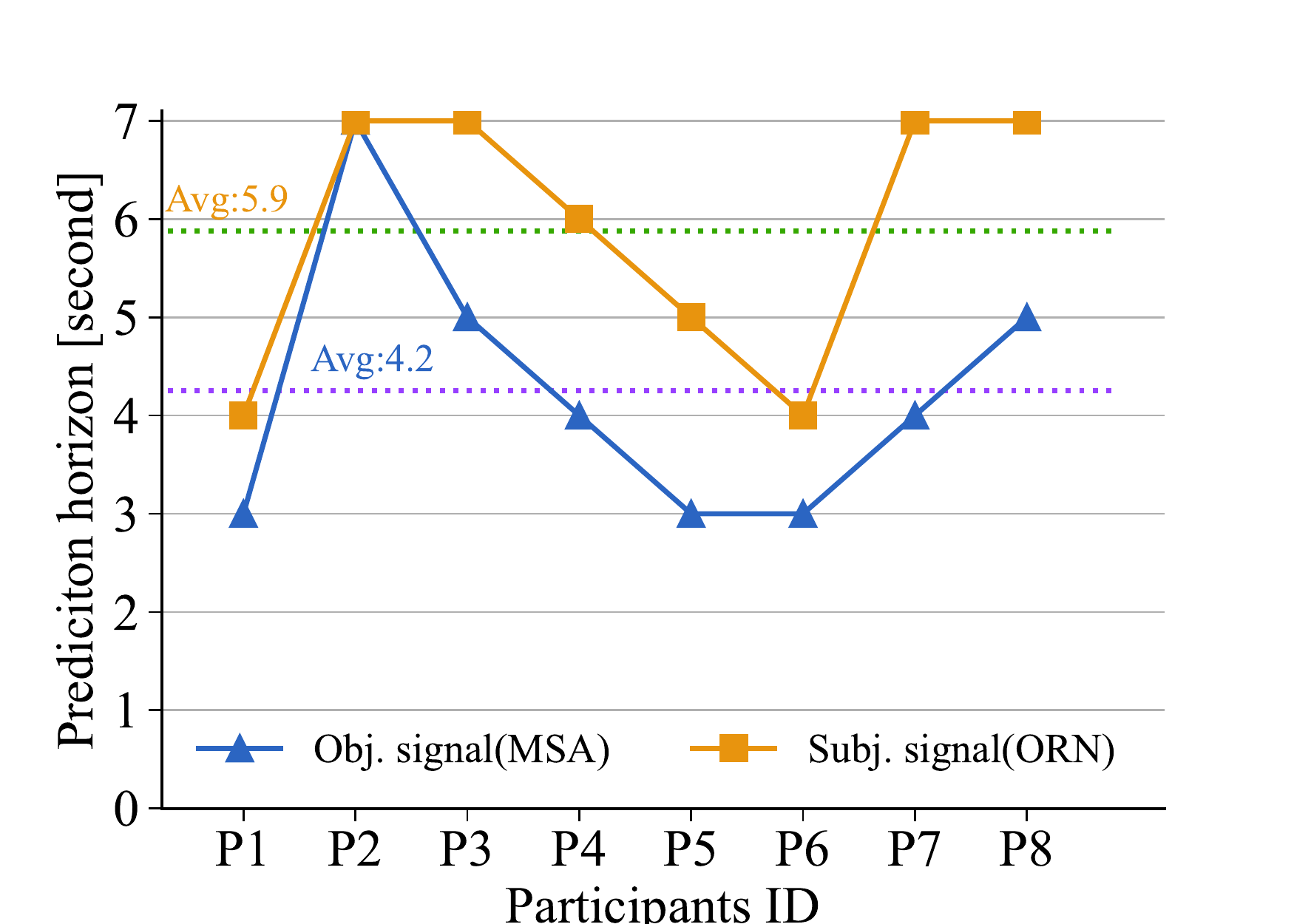}
		\caption{Prediction horizons ($T$)}
		\label{fig_res_pred_hori}
	\end{minipage}	
	\begin{minipage}[t]{0.48\textwidth}
		\centering
		\includegraphics[width=1\textwidth]{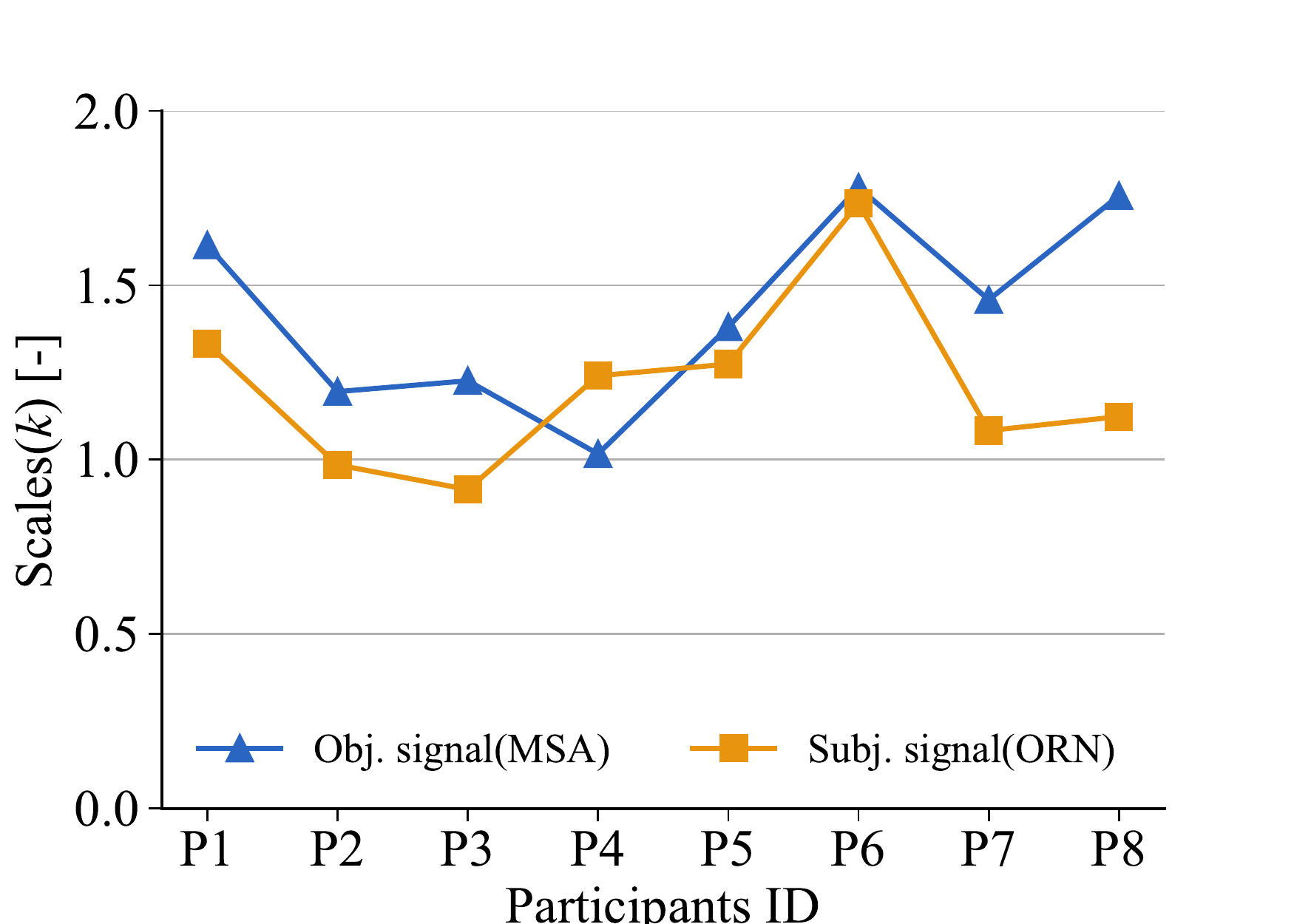}
		\caption{Damage estimation coefficient ($k$)}
		\label{fig_res_damage_est}
	\end{minipage}
\end{figure}

The damage estimation coefficients ($k$) on objective signals showed differences among drivers as shown in Figure \ref{fig_res_damage_est}. It's natural that the tendency of $k$ on objective signals was consistent with the variations of maximum MSA in Figure \ref{fig_data_preview}. Participant P4 had the lowest damage estimation, i.e., lowest sensitivity to collision damage, while participant P8 would put more caution on possible collision harm. The tendency of $k$ on subjective signals can not support the damage estimation analysis as had been illustrated above. However, its tendency was similar to objective signals, implying that the coherence between subjective perception and objective actions.

\subsection{Individual Models on Objective signals}
\begin{figure}[htbp]
	\centering
	\subfigure[Fitted PODAR and raw objective signals]{
		\includegraphics[width=0.98\textwidth]{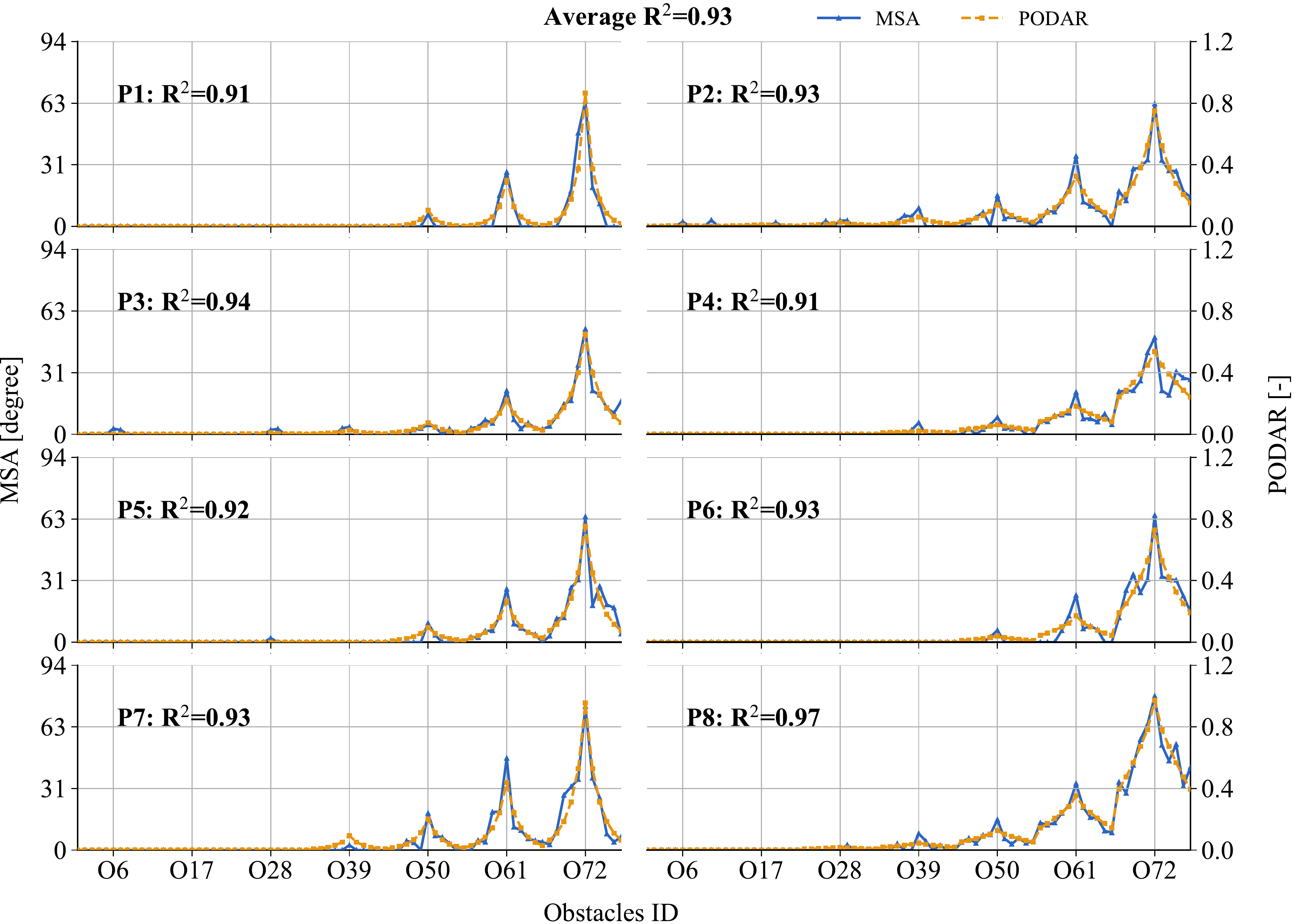}}
	\subfigure[temporal attenuation ($e^{-A\cdot t}$)]{
		\includegraphics[width=0.48\linewidth]{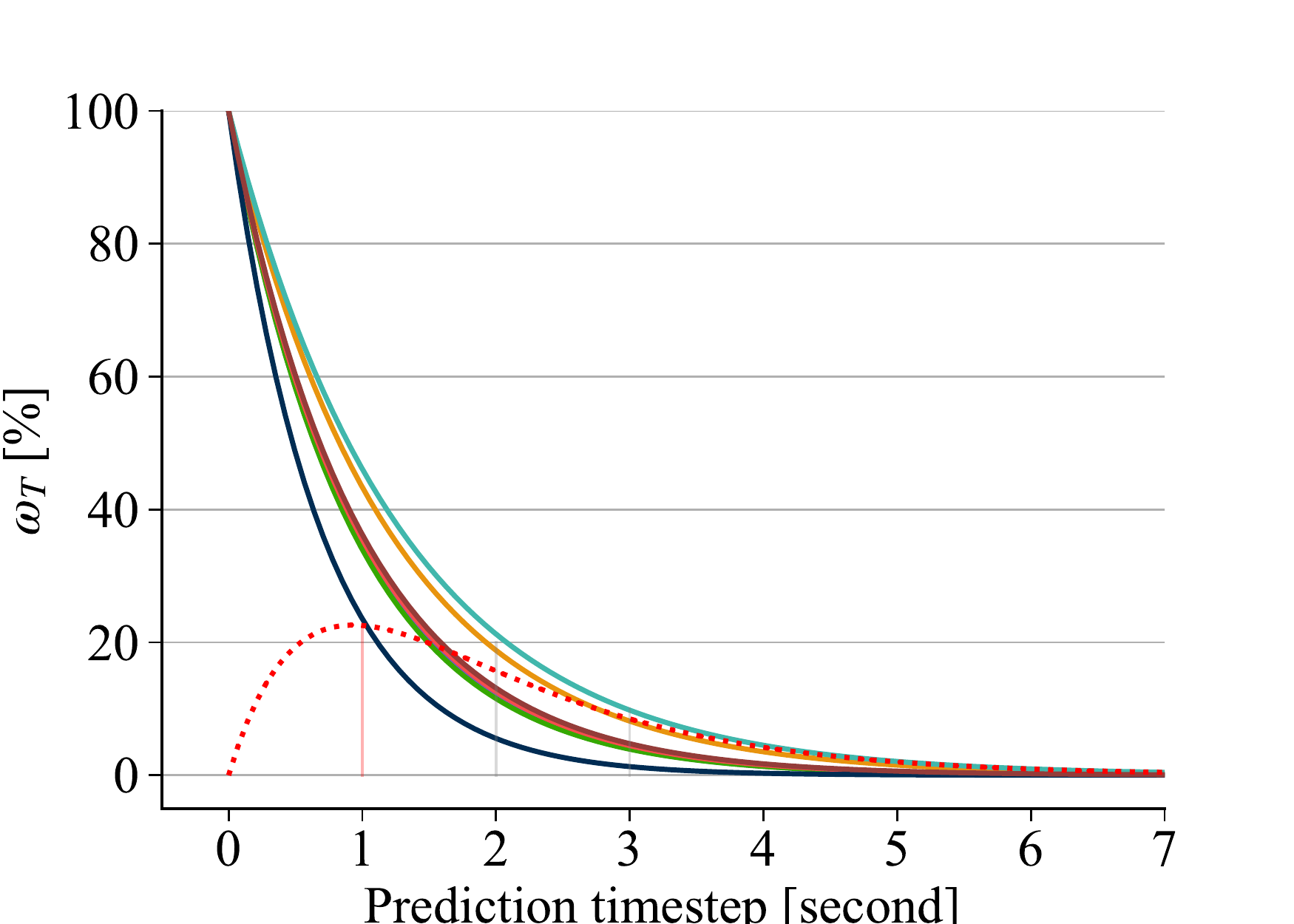}}
	\subfigure[spatial attenuation ($e^{-B\cdot d}$)]{
		\includegraphics[width=0.48\linewidth]{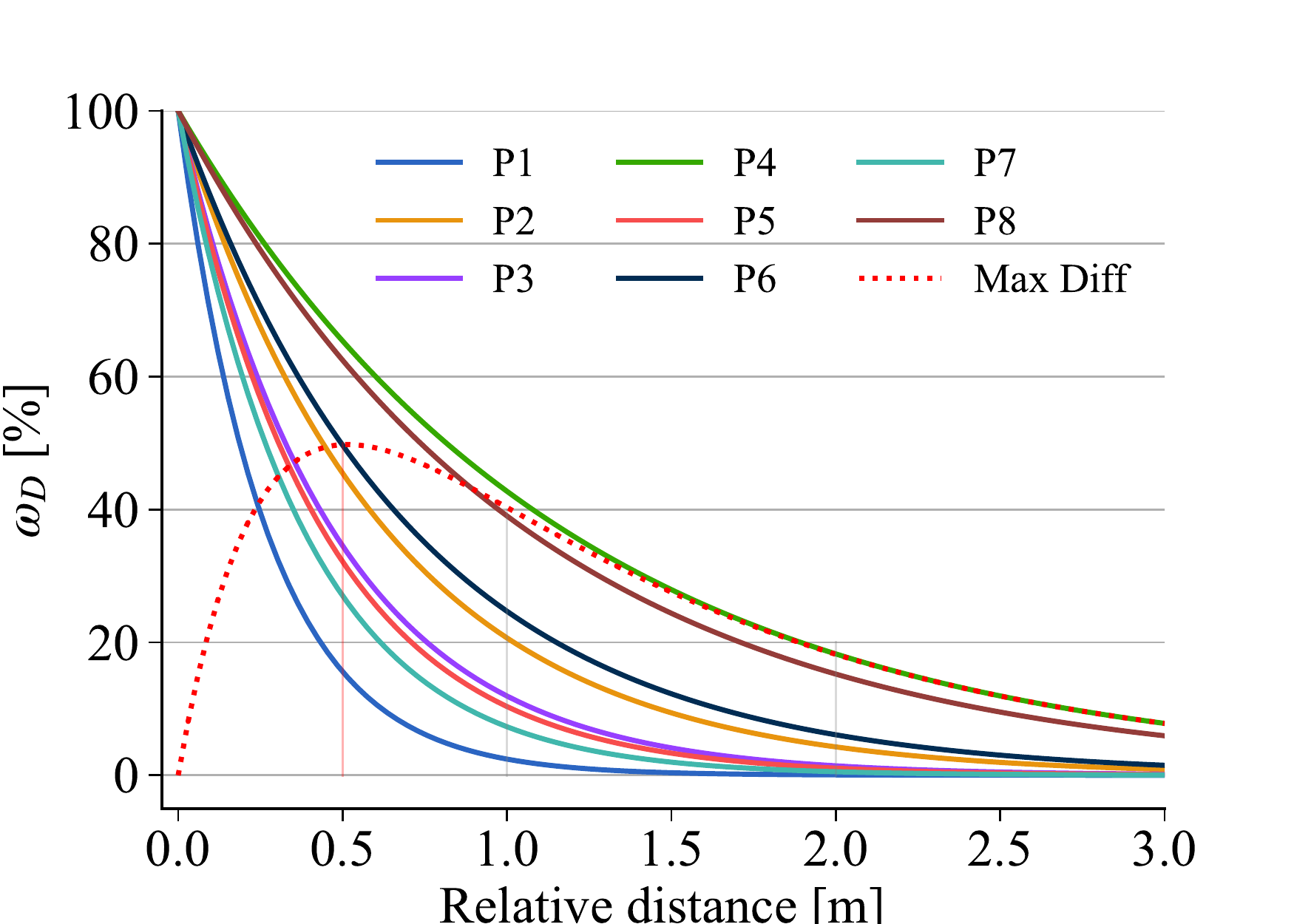}}
	\begin{tabular}{ c c c c c c c c c}
		\hline 
          & P1 & P2 & P3 & P4 & P5 & P6 & P7 & P8 \\ 
        \hline 
        Estimated $A$ & 1.060 & 0.836 & 1.076 & 1.079 & 1.036 & 1.446 & 0.774 & 1.015 \\ 
        \hline 
        Estimated $B$ & 3.716 & 1.577 & 2.126 & 0.850 & 2.269 & 1.399 & 2.617 & 0.941 \\ 
		\hline 
	\end{tabular} 
	\caption{Individual PODAR model of each participant on objective signals}
	\label{fig_res_obj}
\end{figure}

The outputted PODAR values from each individual model and raw MSA values corresponding to the 77 obstacles, as well as the $\text{R}^2$ indexes, are shown in Figure \ref{fig_res_obj}a. Figure \ref{fig_res_obj}b and \ref{fig_res_obj}c show the temporal and spatial attenuation functions ($\omega_T$ and $\omega_D$) of the eight drivers according to calibrated parameters $A$ and $B$. The max differences between the drivers are also calculated and drawn as the red dash line. The values of $A$ and $B$ are also provided in the table of Figure \ref{fig_res_obj}.

The PODAR models perform well on the fitness of raw MSA values with a minimum $\text{R}^2=0.91$ and average $\text{R}^2=0.93$. The assumption of symmetry in the tendency of the risk from sub-sequence obstacles was verified, which indicates that driving risk perception is independent of obstacle direction on a straight road. 

Participant P6 had the steepest descent curve (the biggest $A$ value) in the temporal dimensionality, and also can be observed by examining the PODAR tendency of the main-sequence obstacles in Figure \ref{fig_res_obj}a. Participant P6 would be concerned about the oncoming obstacles in the near future but not be sensitive to the far future damage. On the contrary, except for the near future collision, participants P2 and P7 also kept more attention to the far future possible damages. Generally, the attention on the future damage will decrease to less than 20\% after 2 seconds and less than 10\% after 3 seconds for nearly all the participants. The maximum difference is at around predictive 1 second, where the difference in the risk estimates of different drivers for this period is about 20\%.

Participant P1 had the least sensitive (the biggest $B$ value) to driving risk perception in the spatial dimensionality, which can also be observed in the PODAR tendency of the sub-sequence obstacles in Figure \ref{fig_res_obj}a. That means that participant P1 tended to maintain a small distance to the obstacles, while Participants P4 and P8 would drive far away (smaller $B$ values) from the lateral obstacles to keep safe. Generally, the attention to the damages from far-way obstacles will decrease to less than 40\% and 20\% once the relative distance is larger than 1m and 2m, respectively. Extremely, participant P1 will even be indifferent to objects whose relative distance to the host vehicle is farther than 1m. The attenuation difference between participants in spatial dimensionality varied more than in temporal dimensionality. The perceived risk differences to obstacles at 0.2-1m are over 40\% for different drivers, which can infer a hypothesis that the treatment for obstacles at 0.2-1m will have huge effect on the drivers' acceptance or trust of autonomous vehicles.

\subsection{Subjective Indicator}
\begin{figure}
	\centering
	\subfigure[Fitted PODAR and raw subjective signals]{
		\includegraphics[width=0.98\textwidth]{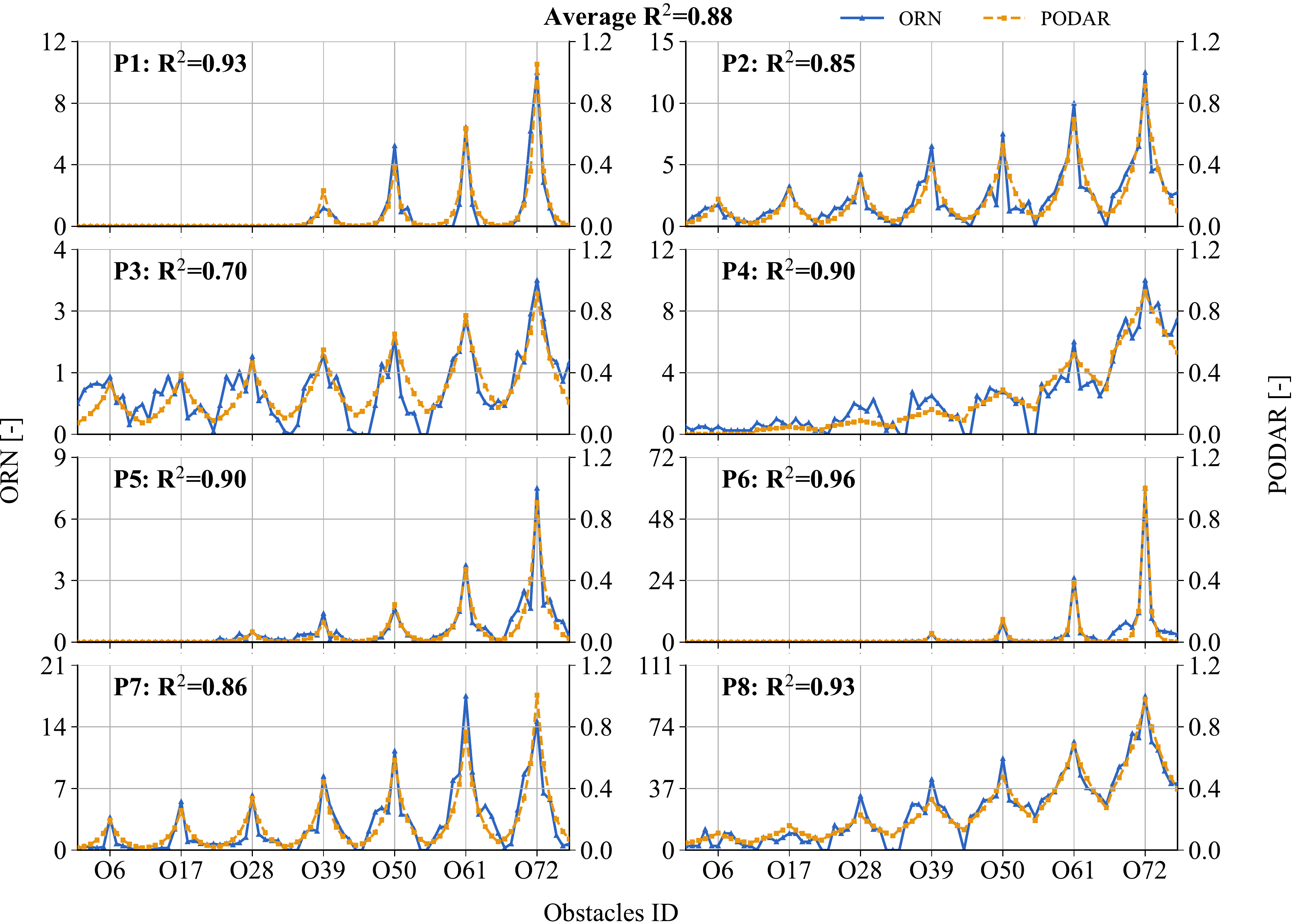}}
	\subfigure[temporal attenuation ($e^{-A\cdot t}$)]{
		\includegraphics[width=0.48\linewidth]{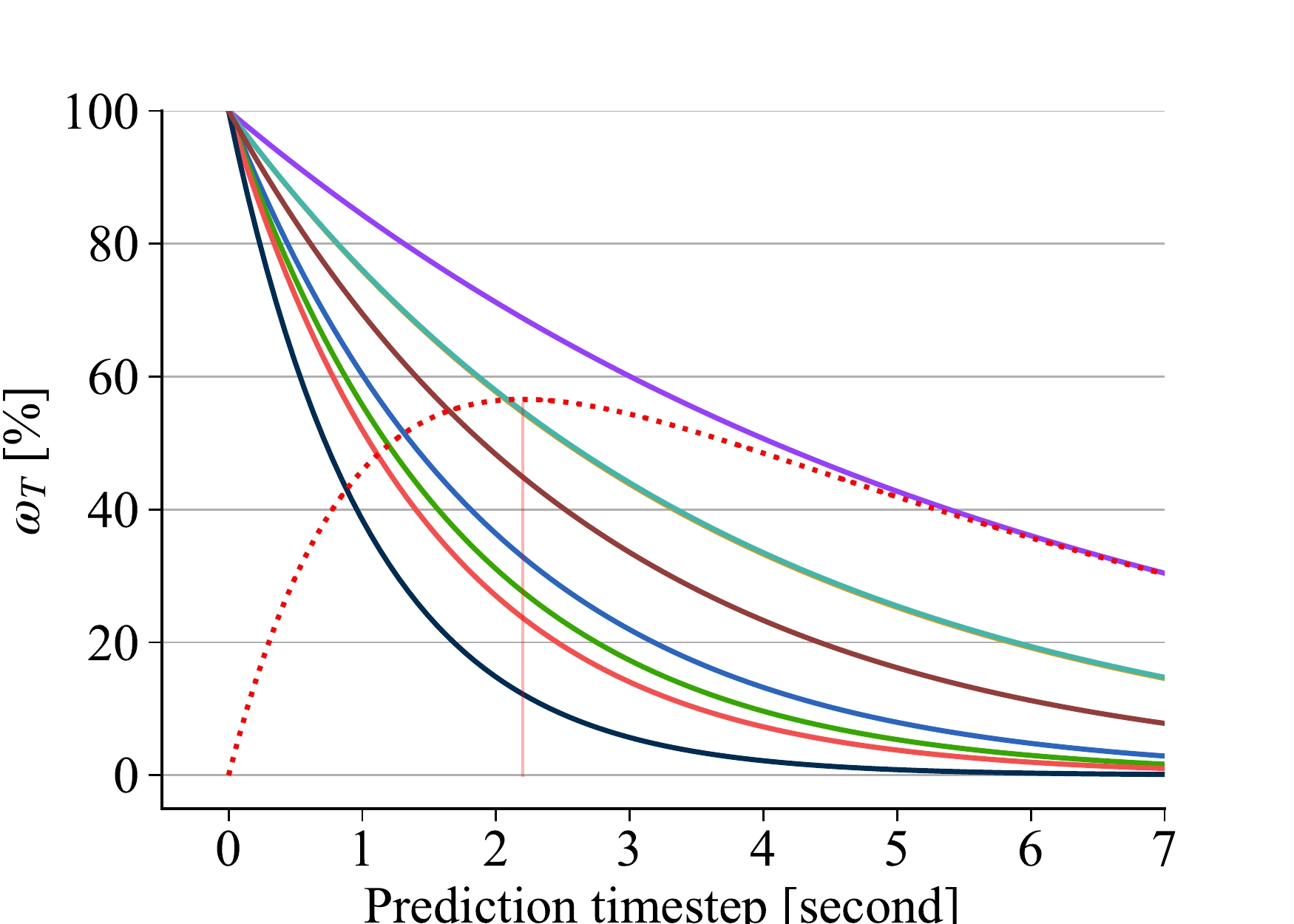}}
	\subfigure[spatial attenuation ($e^{-B\cdot d}$)]{
		\includegraphics[width=0.48\linewidth]{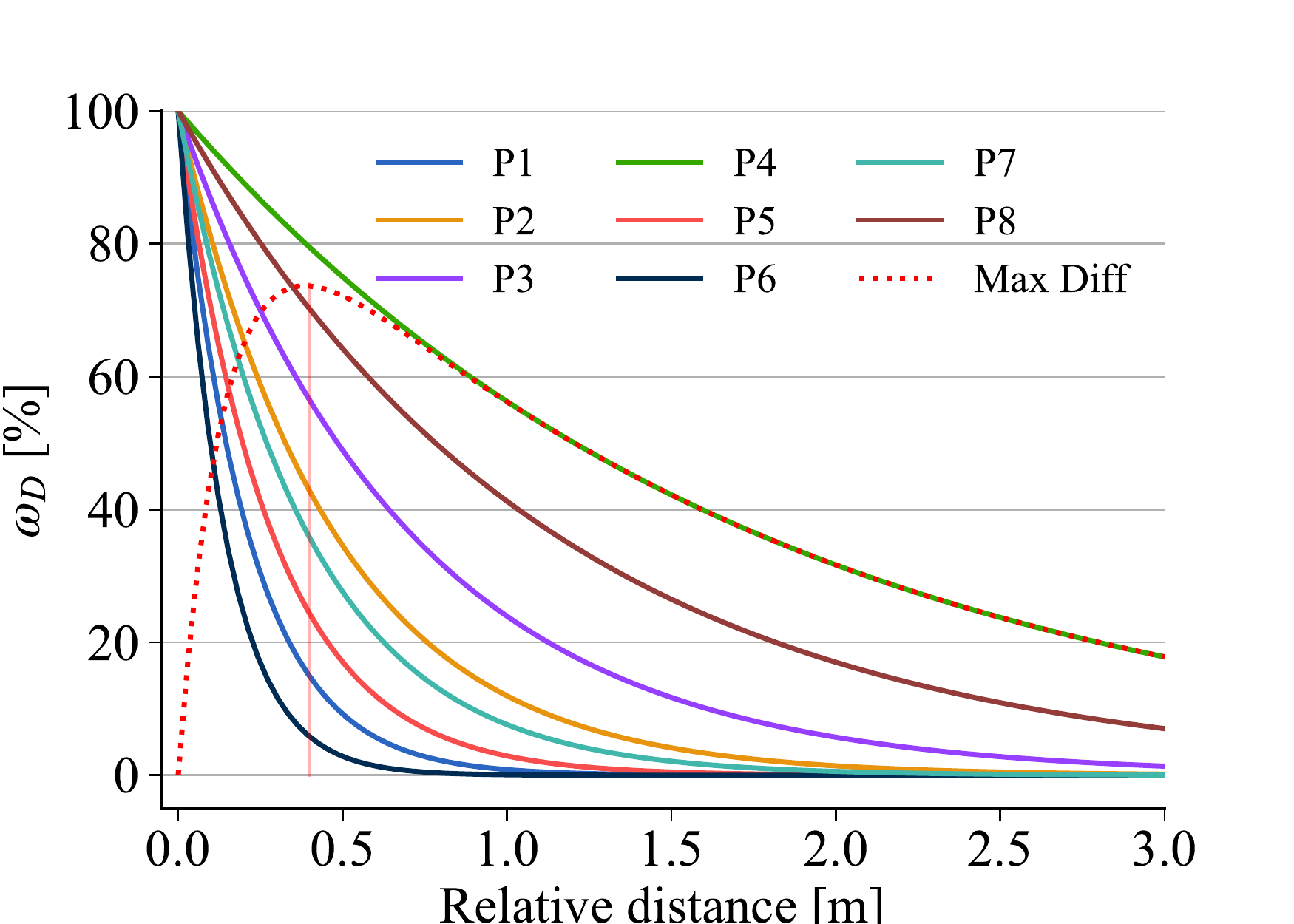}}	
	\begin{tabular}{c c c c c c c c c}
		\hline 
          & P1 & P2 & P3 & P4 & P5 & P6 & P7 & P8 \\ 
        \hline 
        Estimated $A$ & 0.506 & 0.275 & 0.170 & 0.585 & 0.655 & 0.956 & 0.273 & 0.364 \\ 
        \hline 
        Estimated $B$ & 4.765 & 2.125 & 1.430 & 0.575 & 3.537 & 7.115 & 2.573 & 0.886 \\ 
		\hline 
	\end{tabular} 
	\caption{Individual PODAR model of each participant on subjective signals}
	\label{fig_res_sub}
\end{figure}

Figure \ref{fig_res_sub} shows the results of individual PODAR models on subjective signals. The fitness between PODAR values and raw ORNs performed well with an average $\text{R}^2=0.88$. Model of participant P3 had the worst $\text{R}^2=0.70$ value due to the unexpected feedback numbers responding to the far obstacles.

For the temporal attenuation (Figure \ref{fig_res_sub}b), Participant P3 had the smallest $A$ value and would be cautious about the far future collision damages. Participant P3 will not decrease his attention on damages from far future collisions and even keep 40\% damage perception to the collisions in the future 5 seconds. At the same time, participant P6 was insensitive to the future collisions subjectively, like his performance on the objective signals.

Figure \ref{fig_res_sub}c shows the estimated parameters and $B$ in the spatial attenuation function $\omega_D$. Participants P4 and P8 would tend to keep far lateral distances (smaller $B$ values, i.e., keeping high-risk perception about far distance obstacles) to the obstacles subjectively, and they actually acted upon, thus resulting in the smaller $B$ values in the objective signal models (Figure \ref{fig_res_obj}c). Whereas, Participant P6 also acts directly (Figure \ref{fig_res_obj}c) even he/she subjectively perceived very low risks to obstacles at 0.5m. 

Compared with objective signals, the differences between drivers in subjective signals are more huge. The maximum differences for subjective signals are nearly 60\% in temporal dimensionality and 70\% in spatial dimensionality, while 20\% and 50\% for objective signals. Besides, the differences in temporal dimensionality are at predictive 1 second for objective signals, while at predictive 2 seconds for subjective signals.

\subsection{Comparison between objective and subjective model parameters}
\begin{figure}[t]
	\centering
	\subfigure[parameter $A$]{
		\includegraphics[width=0.48\linewidth]{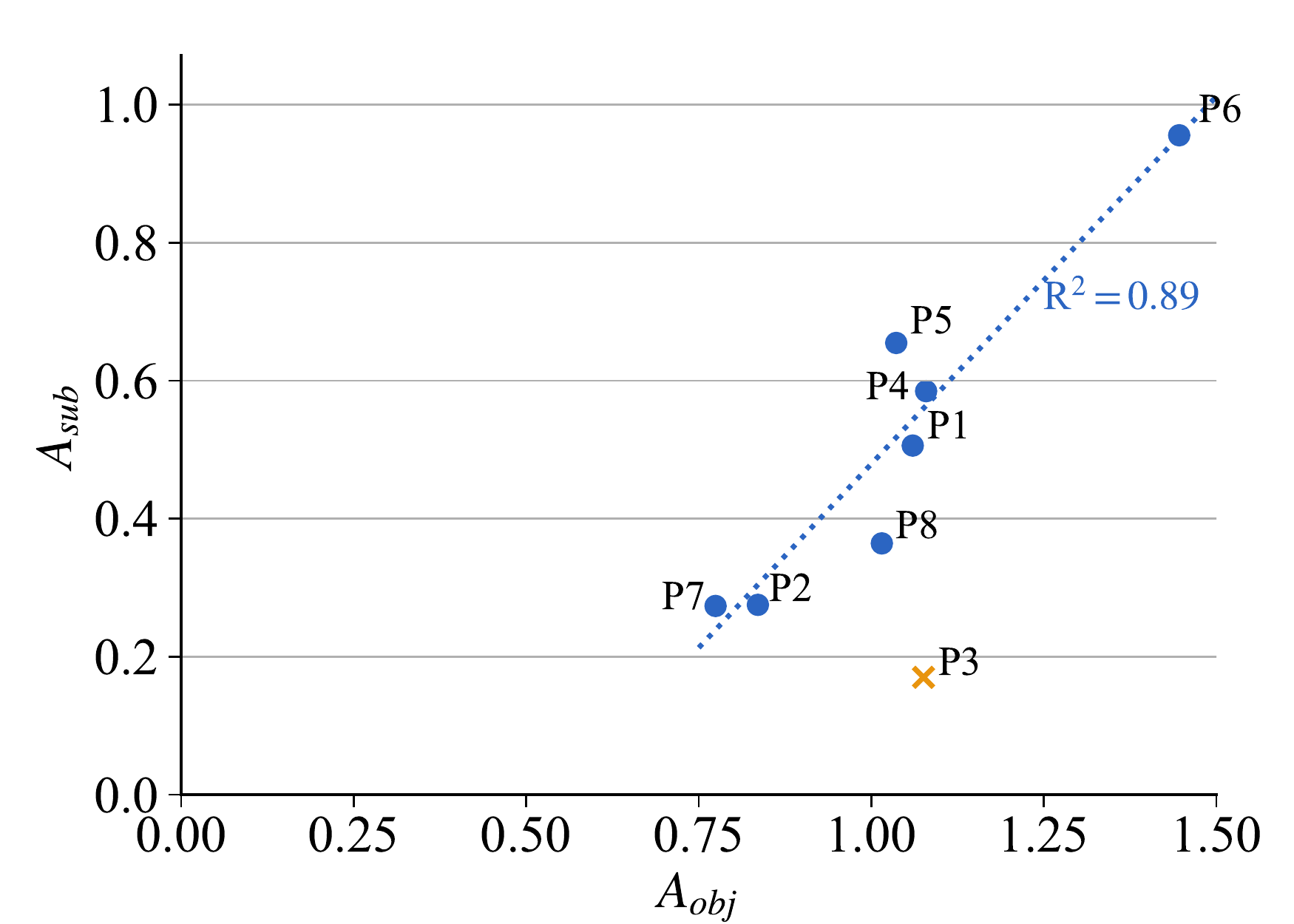}}
	\subfigure[parameter $B$]{
		\includegraphics[width=0.48\linewidth]{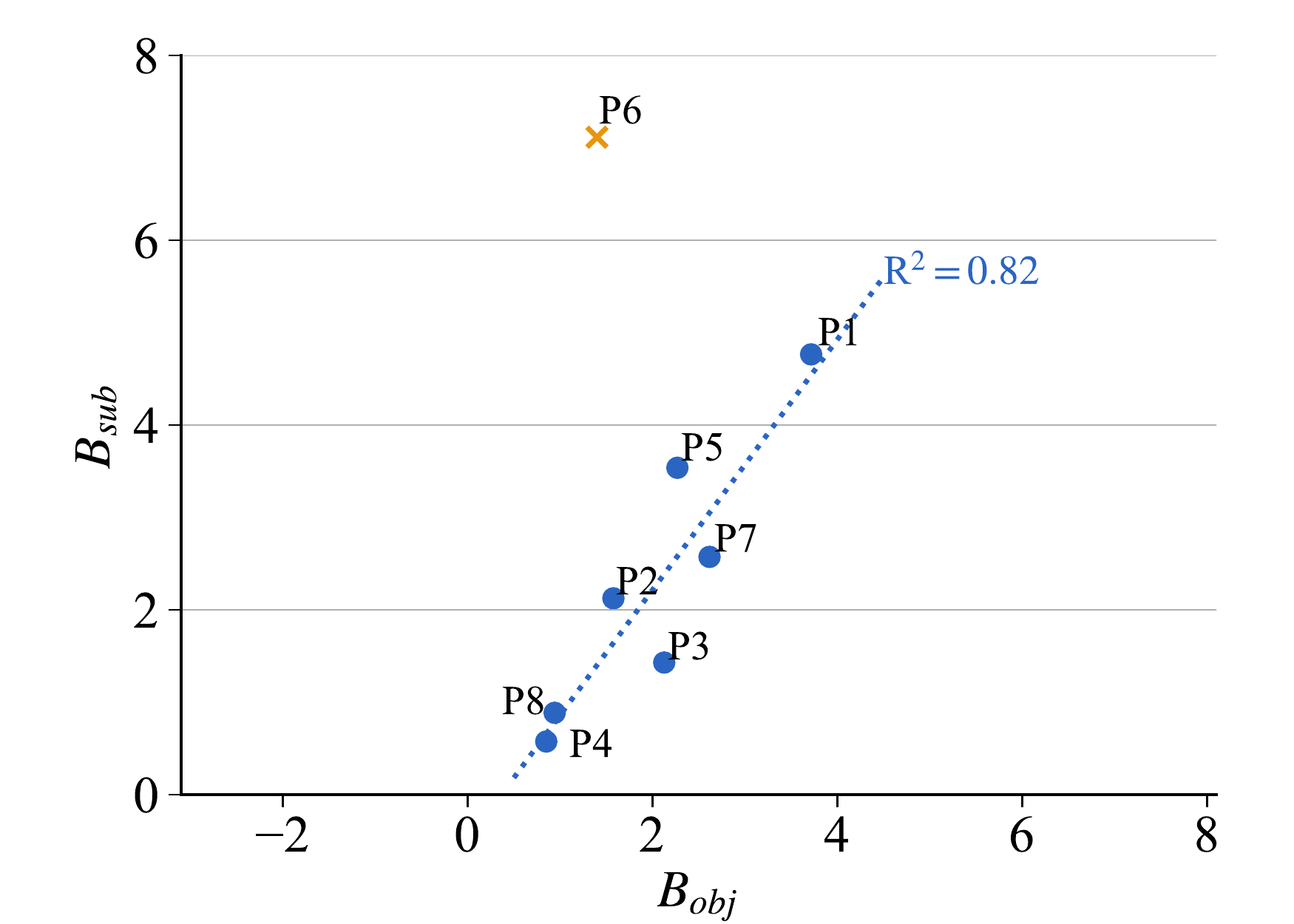}}
	\caption{Comparison of parameters $A$ and $B$ between objective and subjective signals}
	\label{fig_AB_vs}
\end{figure}

For each participant, there was a strong association between $A$ and $B$ values from their objective and subjective PODAR models, as shown in Figure \ref{fig_AB_vs}. Also, there were unexpected outlier points, participant P3 for parameter $A$ and P6 for parameter $B$. Generally, the subjectively perceived risk would induce corresponding actions. However, there were some differences in the temporal and spatial dimensionalities. The $A$ values from objective PODAR models were bigger than that from subjective PODAR values. Recalling that smaller values mean a more slowly declining of the damage attenuation curves, these results indicate that although drivers perceive far future risks, they do not instantly react at the current moment. The parameter $B$, in turn, had similar numerical values between objective and subjective signals, which means that drivers would act directly to deal with risk from spatial dimensionality.

\section{Discussion}
\label{discussion}
The contribution of this paper is to verify the feasibility of the PODAR model in the discrimination of individual characteristics on driving risk perception among different drivers, which provides a concise yet effective way to understand how drivers assess driving risks and why there are differences. It only takes four parameters to build a customized risk perception model, providing a scalar risk value and suitable for general scenes by resolving temporal and spatial dimensionalities.

Open-accessed data from a previous study\cite{ref3}, which were used to verify the proposed Driver's Risk Field (DRF) model, is used in our research.  Unlike PODAR, the DRF model used a field-based method to separate the room in front of the host vehicle into lateral and longitudinal directions and utilized the Gaussian function (two parameters needed to be calibrated) and power function (three parameters needed to be calibrated) to fit respectively. However, signals from obstacles in each column or each row in Figure \ref{fig_DS_design}a were fitted, leading to 7 Gaussian and 11 power functions for each driver. Totally $2\times7+3\times11=47$ parameters were needed in the modelling for each participant, as shown in Table \ref{tab_comparison}. The DRF model fitness indexes $\text{R}^2$ of objective and subjective signals were 0.77 and 0.69 in the lateral direction and 0.86 and 0.98 in the longitudinal direction. Compared with the calibration results of the PODAR model ($\text{R}^2$ are 0.93 and 0.88 for objective and subjective signals respectively), DRF only is superior to the fitness of subjective signals in the longitudinal direction. However, it should be noted that the $\text{R}^2$ of the DRF model was the average of several, meaning it used parameters (47) are nearly 12 times that (4) of the PODAR model. Therefore, it was expected that the fitness of the DRF model should be better than PODAR, while the truth is that PODAR is better in most aspects with fewer parameters.

The favorable performance of the PODAR model benefits from its deep insight of the nature of driving risk. As mentioned in our previous work\cite{ref_podar}, the risk is regarded as the projection of collision damage, which numerically encourages the PODAR to find a function to describe the attenuation from spatial and temporal dimensionalities. The exponential function with $e$ as base came into our sight because $e$ is the natural constant and the exponential function naturally meets the requirements of Equation \ref{eq_attenuation}. Compared with the Gaussian and power function used in the DRF model, the results show the better capacity of PODAR model in risk estimation. Besides, the way building the PODAR model also assigns its parameters,  $A$, $B$, $T$, and $k$, with actual physical meanings. In contrast, parameters used in the DRF model only support a qualitative comparison instead of a quantitative analysis. The interpretable method allows establishing different risk perception models only by simply modifying the parameters, on the one hand, providing a way to grasp the differences in risk perception among drivers, on the other hand, supporting the automated vehicles to generate individual customized behaviors by introducing driver's PODAR model into the decision and planning algorithms optimization process.

\begin{table}[b]
	\centering
	\caption{Performance comparison between PODAR and DRF}
	\label{tab_comparison}
	\begin{tabular}{c c c c c c c}
		\hline
		Model & Signals & Directions & \tabincell{c}{Num. of \\ fitted data} & $\text{R}^2$ &  \tabincell{c}{Num. of parameters \\ in each function} & \tabincell{c}{Total parameters\\for each driver} \\ 
		\hline
		\multirow{4}*{DRF} & \multirow{2}*{objective} & longitudinal & 7 & 0.98 & 3 & \multirow{2}*{47} \\
		\cline{3-6}
						~  & ~ & lateral & 11 & 0.69 & 2 & ~\\
		\cline{2-7}
						~  & \multirow{2}*{subjective} & longitudinal & 7 & 0.86 & 3 & \multirow{2}*{47}\\
		\cline{3-6}
						~  & ~ & lateral & 11 & 0.77 & 2 & ~\\
		\hline
		\multirow{2}*{PODAR} & \multicolumn{2}{c}{objective} & 77 & 0.93 & 4 & 4 \\
		\cline{2-7}
						   ~ & \multicolumn{2}{c}{subjective} & 77 & 0.88 & 4 & 4 \\
		\hline
	\end{tabular}
\end{table}

Due to the limitation of the available dataset, this paper only makes a calibration in the static obstacles scenarios. However, it has shown massive potential in modelling human-like driving risk perception with features that are interpretable, concise, and practical. One key class of scenes, i.e., the lateral conflict condition like collisions in a 90-degree intersection, although have been proved to can be modeled by PODAR with a numerical simulation, we will further take calibration in the future work.

From a conceptual point of view, the core of PODAR model is the basic framework as shown in Equation \ref{eq_PODAR_struc}, \ref{eq_damage} and \ref{eq_attenuation}. Using the exponential function is a good choice, and temporal parameter $A$ and spatial parameter $B$ were independent and calibrated separately. However, we do not exclude the methods coupling parameters $A$ and $B$ to get a more accurate risk mode as long as it obeys the basic framework of PODAR. Besides, trajectory prediction is also significant in the PODAR modelling, especially when facing complex road alignments and considering the uncertainty of traffic participant behaviors.

The differences in presented risk perception among participants can be described by the PODAR model, while the internal motivation behind people's different cognition still needs more research. In the experiment where the dataset was collected, the obstacle appeared suddenly and brought stimulation to drivers, resulting in exaggerated reactions. This is acceptable in laboratory research to generate valid data, whereas surrounding objects' motions are continuous. Future works will also focus on validating the PODAR model in naturalistic driving situations.

\section{Conclusion}
\label{conclusion}
This paper verified the PODAR model using an open-accessed dataset from an obstacle avoidance experiment on a straight road in a driving simulation. The results proved the potential of the PODAR model in the individual modelling of drivers' risk perception, providing a concise but effective way to explain how drivers perceive and assess risks from road objects. Besides, we also can conclude the following:
\begin{itemize}
	\item The PODAR model can fit the drivers' risk perception well with average $\text{R}^2=0.93$ for the objective signals and average $\text{R}^2=0.88$ for subjective signals, better than the DRF model in previous literature.
	\item Only four physical-meaningful and interpretable parameters, $T$, $k$, $A$, and $B$, corresponding to the prediction horizon, potential damage scale, temporal attenuation, and spatial attenuation, are needed to model different drivers' risk perception characteristics.
	\item Drivers would not directly react to the risk from a future collision but would act directly to deal with risks caused by the distance factor.
	\item The average prediction horizon is 6 seconds subjectively but is 4 seconds objectively. 
	\item Most drivers would objectively pay less than 20\% risk attention to future collisions in predictive 2 seconds and the obstacles farther than 2m.
	\item The treatment for obstacles at 0.2-1m would significantly affect the drivers' acceptance or trust of autonomous vehicles.
\end{itemize}

\section{Acknowledgement}
This study is supported in part by the National Natural Science Foundation of China under No. 52102411 and No. U20A20334, and in part by the Tsinghua University-Didi Joint Research Center for Future Mobility.



\bibliographystyle{elsarticle-num} 
\bibliography{ref}

\end{sloppypar}
\end{document}